\documentclass[final]{siamltex}
\usepackage{graphicx,array,amssymb,psfrag}

\newcommand{\BEAS}{\begin{eqnarray*}}
\newcommand{\EEAS}{\end{eqnarray*}}
\newcommand{\BEA}{\begin{eqnarray}}
\newcommand{\EEA}{\end{eqnarray}}
\newcommand{\BEQ}{\begin{equation}}
\newcommand{\EEQ}{\end{equation}}
\newcommand{\BIT}{\begin{itemize}}
\newcommand{\EIT}{\end{itemize}}
\newcommand{\BNUM}{\begin{enumerate}}
\newcommand{\ENUM}{\end{enumerate}}

\newcommand{\BA}{\begin{array}}
\newcommand{\EA}{\end{array}}

\newcommand{\ones}{\mathbf 1}

\newcommand{\reals}{{\mbox{\bf R}}}

\newcommand{\symm}{{\mbox{\bf S}}}  

\newcommand{\Card}{\mathop{\bf Card}}
\newcommand{\Tr}{\mathop{\bf Tr}}

\newcommand{\argmin}{\mathop{\rm argmin}}

\begin{document}

\title{Smooth Optimization \\with Approximate Gradient}

\author{Alexandre d'Aspremont\thanks{ORFE Department, Princeton University, Princeton, NJ 08544. \texttt{aspremon@princeton.edu}}}

\maketitle

\begin{abstract}
We show that the optimal complexity of Nesterov's smooth first-order optimization algorithm is preserved when the gradient is only computed up to a small, uniformly bounded error. In applications of this method to semidefinite programs, this means in some instances computing only a few leading eigenvalues of the current iterate instead of a full matrix exponential, which significantly reduces the method's computational cost. This also allows sparse problems to be solved efficiently using sparse maximum eigenvalue packages.
\end{abstract}

\section{Introduction}
In \cite{Nest83} it was shown that smooth convex minimization problems of the form:
\[
\BA{ll}
\mbox{minimize } & f(x)\\
\mbox{subject to} & x \in Q,\\
\EA\]
where $f$ is a convex function with Lipschitz continuous gradient and $Q$ is a sufficiently simple compact convex set, could be solved with a complexity of $O(1/\sqrt{\epsilon})$, where $\epsilon$ is the precision target. Furthermore, it can be shown that this complexity bound is optimal for that class of problems (see \cite{Nest03a} for a dicussion). More recently, \cite{Nest03} showed that this method could be combined with a smoothing argument to produce an $O(1/\epsilon)$ complexity bound for non-smooth problems where the objective has a saddle-function format. In particular, this meant that a broad class of semidefinite optimization problems could be solved with significantly lower memory requirements than interior point methods and a better complexity bound than classic first order methods (bundle, subgradient, etc).

Here, we show that substituting an approximate gradient, which may allow significant computation and storage savings, does not affect the optimal complexity of the algorithm in \cite{Nest83}. It is somewhat intuitive that an algorithm which exhibits good numerical performance in practice should be robust to at least some numerical error in the objective function and gradient computations since all implementations are necessarily computing these quantities up to some multiple of machine precision. Our objective here is to make that robustness explicit in order to design optimal schemes using only approximate gradient information.

For non-smooth problems, when the objective function $f(x)$ can be expressed as a saddle function on a compact set, the method in \cite{Nest03} starts by computing a smooth (i.e. with Lipschitz continuous gradient), uniform $\epsilon$-approximation of the objective function $f(x)$, it then uses the smooth minimization algorithm in \cite{Nest83} to solve the approximate problem. When this smoothing technique is applied to semidefinite optimization, computing exact gradients requires forming a matrix exponential, which is often the dominant numerical step in the algorithm. 

Although there are many different methods for computing this matrix exponential (see \cite{Mole03} for a survey), their complexity is comparable to that of a full eigenvalue decomposition of the matrix. In problem instances where only a few leading eigenvalues suffice to approximate this exponential, the per iteration complexity of the algorithm described here becomes comparable to that of classical first-order methods such as the bundle method (see \cite{helm00}) or subgradient methods (see \cite{Shor85} for example), which have a global complexity bound of $O(1/\epsilon^2)$ (see \cite{Nest03a}), while keeping the optimal complexity of $O(1/\epsilon)$ of the algorithm in \cite{Nest03}. 

We apply this result to a maximum eigenvalue minimization problem (or semidefinite program with constant trace). We first recall the complexity bound derived in \cite{Nest04a} based on a smoothing argument, using exact gradients. We produce a rough theoretical estimate of the number of eigenvalues required for convergence when approximate gradients are used. We then derive an explicit condition on the quality of the gradient approximation to guarantee convergence and compute a bound on the number of iterations. We show both on randomly generated problem instances and on problems generated from biological data sets that actual computational savings vary significantly with problem structure but can be substantial in some cases.

The paper is organized as follows. In the next section, we prove convergence of the algorithm in \cite{Nest83} when only an approximate gradient is used. In Section \ref{s:semidef} we describe how these results can be applied to semidefinite optimization. Finally, in the last section we test their performance on semidefinite relaxation and maximum eigenvalue minimization problems.

\section{Smooth optimization with approximate gradient}
\label{s:algo} Following the results and notations in \cite[\S 3]{Nest03}, we study the problem:
\BEQ
\label{eq:smooth-min}
\BA{ll}
\mbox{minimize } & f(x)\\
\mbox{subject to} & x \in Q,\\
\EA\EEQ 
where $Q\subset \reals^n$ is a compact convex set and $f$ is a convex function with Lipschitz continuous gradient, such that:
\[
\|\nabla f(x) - \nabla f(y)\|^* \leq L \|x-y\|,\quad x,y\in Q,
\]
for some $L>0$, which also means:
\BEQ
\label{eq:lip-ineq}
f(y) \leq f(x) +\langle \nabla f(x),y-x\rangle +\frac{1}{2}L\|y-x\|^2, \quad x,y\in Q.
\EEQ
The key difference here is that the \emph{oracle} information we obtain for $\nabla f$ is \emph{noisy}. Note that the function values are not required to compute iterates in the algorithm described here, so even if our knowledge of function values $f(x)$ is noisy, we will always use  exact values in the proofs that follow. At each iteration, we obtain $\tilde \nabla f(x)$ satisfying:
\BEQ
\label{eq:approx-grad}
|\langle \tilde \nabla f(x)-\nabla f(x),y-z\rangle| \leq \delta \quad x,y,z\in Q,
\EEQ
for some precision level $\delta >0$. Throughout the paper, we assume that $Q$ is simple enough so that this condition can be checked efficiently. As in \cite{Nest83}, we also assume that certain projection operators on $Q$ can be computed efficiently and we refer the reader to the end of this section for more details. Here, $d(x)$ is a prox-function for the set $Q$, i.e. continuous and strongly convex on $Q$ with parameter $\sigma$ (see \cite{Nest03a} or \cite{Hiri96} for a discussion of regularization techniques using strongly convex functions). We let $x_0$ be the center of $Q$ for the prox-function $d(x)$ so that:
\[
x_0\triangleq\argmin_{x\in Q}d(x),
\]
assuming w.l.o.g. that $d(x_0)=0$, we then have:
\BEQ
\label{eq:d-strong-convex}
d(x)\geq\frac{1}{2}\sigma\|x-x_0\|^2.
\EEQ
We denote by $\tilde T_Q(x)$ a solution to the following subproblem:
\BEQ
\label{eq:tq-tilde}
\tilde T_Q(x) \triangleq \argmin_{y\in Q}\left\{\langle \tilde \nabla f(x),y-x \rangle+\frac{1}{2}L\|y-x\|^2\right\}
\EEQ
We let $y_0=\tilde T_Q(x_0)$ where $x_0$ is defined above. We recursively define three sequences of points: the current iterate $\{x_k\}$, the corresponding $y_k=\tilde T_Q(x_k)$, together with
\BEQ
\label{eq:zk-def}
z_k \triangleq \argmin_{x \in Q}\left\{\frac{L}{\sigma}d(x)+\sum_{i=0}^k\alpha_i[ f(x_i)+\langle \tilde \nabla f(x_i),x-x_i\rangle]\right\}
\EEQ
and a step size sequence $\{\alpha_k\}\geq 0$ with $\alpha_0\in(0,1]$ so that
\BEQ
\label{eq:xyz-update}
\BA{l}
x_{k+1}=\tau_k z_k+(1-\tau_k)y_k\\
y_{k+1}=\tilde T_Q(x_{k+1})\\
\EA
\EEQ
where $\tau_k=\alpha_{k+1}/A_{k+1}$ with $A_k=\sum_{i=0}^k\alpha_i$. We implicitly assume here that the two subproblems defining $y_k$ and $z_k$ can be solved very efficiently (in the examples that follow, they amount to Euclidean projections). We will show recursively that for a good choice of step sequence $\alpha_k$, the iterates $x_k$ and $y_k$ satisfy the following relationship (denoted by $\mathcal{R}_k$):
\[
A_k f(y_k)\leq\psi_k+ A_k g(k,\delta) \quad (\mathcal{R}_k)
\]
where $g(k,\delta)$ measures the accumulated gradient approximation error and will be bounded in Lemma \ref{th:psik}, and
\[
\psi_k\triangleq \min_{x \in Q}\left\{\frac{L}{\sigma}d(x)+\sum_{i=0}^k\alpha_i[ f(x_i)+\langle \tilde \nabla f(x_i),x-x_i\rangle]\right\}.
\]
First, using $d(x)\geq\frac{1}{2}\sigma\|x-x_0\|^2$, then inequality (\ref{eq:lip-ineq}) and condition (\ref{eq:approx-grad}), we have:
\[
\psi_0=\min_{x \in Q}\left\{\frac{L}{\sigma}d(x)+\alpha_0[ f(x_0)+\langle \tilde \nabla f(x_0),x-x_0\rangle]\right\}
\geq \alpha_0 f(y_0) -  \alpha_0 \delta
\]
which is $\mathcal{R}_0$. We can then bound the approximation error in the following result.
\begin{lemma}
\label{th:psik}
Let $\alpha_k$ be a step sequence satisfying:
\BEQ\label{eq:alphak}
0<\alpha_0\leq1 \quad \mbox{and} \quad \alpha_k^2\leq A_k,\quad k\geq 0,
\EEQ
suppose that $(\mathcal{R}_k)$ holds with $x_{k+1}$ and $y_{k+1}$ are defined as in (\ref{eq:xyz-update}), then $(\mathcal{R}_{k+1})$ holds with:
\[
g(k+1,\delta)=(1-\tau_k)g(k,\delta)+\tau_k 3\delta,
\]
where $\tau_k\in[0,1]$ and $g(0,\delta)=\alpha_0\delta$.
\end{lemma}
\begin{proof}
Let us assume that $(\mathcal{R}_k)$ holds. Because $d(x)$ is strongly convex with parameter $\sigma$, the function:
\[
\frac{L}{\sigma}d(x)+\sum_{i=0}^{k}\alpha_i[ f(x_i)+\langle \tilde \nabla f(x_i),x-x_i\rangle]
\]
is strongly convex with parameter $L$. Using this property and the definition of $z_k$ we obtain:
\BEAS
\psi_{k+1}& = & \min_{x \in Q}\left\{\frac{L}{\sigma}d(x)+\sum_{i=0}^{k+1}\alpha_i[ f(x_i)+\langle \tilde \nabla f(x_i),x-x_i\rangle]\right\}\\
& \geq& \min_{x \in Q}\left\{\psi_k+\frac{1}{2}L\|x-z_k\|^2+\alpha_{k+1}[ f(x_{k+1})+\langle \tilde \nabla f(x_{k+1}),x-x_{k+1}\rangle]\right\}.\\
\EEAS
Now, using $(\mathcal{R}_k)$ then the convexity of $f(x)$, we get:
\setlength{\extrarowheight}{1.5ex}
\[\BA{l} 
\psi_k+A_kg(k,\delta)+\alpha_{k+1}[ f(x_{k+1})+\langle \tilde \nabla f(x_{k+1}),x-x_{k+1}\rangle]\\
\geq A_k f(y_k)+\alpha_{k+1}[ f(x_{k+1})+\langle \tilde \nabla f(x_{k+1}),x-x_{k+1}\rangle]\\
\geq A_k [f(x_k)+\langle \nabla f(x_{k+1}),y_k-x_{k+1}\rangle]+\alpha_{k+1}[ f(x_{k+1})+\langle \tilde \nabla f(x_{k+1}),x-x_{k+1}\rangle],\\
\EA\]
and condition (\ref{eq:approx-grad}), together with (\ref{eq:xyz-update}) imply that:
\setlength{\extrarowheight}{1.5ex}
\[\BA{l} 
A_k [f(x_k)+\langle \nabla f(x_{k+1}),y_k-x_{k+1}\rangle]+\alpha_{k+1}[ f(x_{k+1})+\langle \tilde \nabla f(x_{k+1}),x-x_{k+1}\rangle]\\
\geq A_{k+1} f(x_{k+1})+\langle \nabla f(x_{k+1}),A_k y_k- A_k x_{k+1}+\alpha_{k+1}(x-x_{k+1})\rangle-\alpha_{k+1}\delta\\
= A_{k+1} f(x_{k+1})+\alpha_{k+1} \langle \nabla f(x_{k+1}),x-z_{k}\rangle - \alpha_{k+1}\delta.\\
\EA\]
Because $\alpha_k$ satisfies (\ref{eq:alphak}), we have $\tau^2_k\leq A_{k+1}^{-1}$ and can combine the last three inequalities to get:
\BEQ\label{eq:midineq}
\BA{rl}
\psi_{k+1}\geq & A_{k+1}f(x_{k+1})-A_k g(k,\delta)- \alpha_{k+1}\delta\\
& +\min_{x \in Q}\left\{\frac{1}{2}L\|x-z_k\|^2 + \alpha_{k+1} \langle \nabla f(x_{k+1}),x-z_{k}\rangle \right\}\\
\geq & A_{k+1}[f(x_{k+1})-(1-\tau_k) g(k,\delta)- \tau_{k}\delta\\
&+\min_{x \in Q}\left\{\frac{1}{2}L\tau^2_k\|x-z_k\|^2+\tau_{k}\langle \nabla f(x_{k+1}),x-z_{k}\rangle \right\}].\\
\EA\EEQ
Let us define $y \triangleq \tau_k x + (1-\tau_k)y_k$ so that $y-x_{k+1}=\tau_k(x-z_k)$, with:
\[\BA{l}
\min_{x \in Q}\left\{\frac{1}{2}L\tau^2_k\|x-z_k\|^2+\tau_{k}\langle \nabla f(x_{k+1}),x-z_{k}\rangle \right\}\\
=\min_{\{y \in \tau_kQ+(1-\tau_k)y_k\}}\left\{\frac{1}{2}L\|y-x_{k+1}\|^2+\langle \nabla f(x_{k+1}),y-x_{k+1}\rangle \right\}.\\
\EA\]
Combining condition (\ref{eq:approx-grad}) with the fact that $y-x_{k+1}=\tau_k(x-z_k)$ for some $x,z_k\in Q$, we get:
\[\BA{l}
\min_{\{y \in \tau_kQ+(1-\tau_k)y_k\}}\left\{\frac{1}{2}L\|y-x_{k+1}\|^2+\langle \nabla f(x_{k+1}),y-x_{k+1}\rangle \right\}\\
\geq \min_{\{y \in \tau_kQ+(1-\tau_k)y_k\}}\left\{\frac{1}{2}L\|y-x_{k+1}\|^2+\langle \tilde \nabla f(x_{k+1}),y-x_{k+1}\rangle \right\}- \tau_k \delta\\
\EA\]
Now, because $Q$ is convex, we must have $\tau_kQ+(1-\tau_k)y_k \subset Q$ and:
\[\BA{l}
\min_{\{y \in \tau_kQ+(1-\tau_k)y_k\}}\left\{\frac{1}{2}L\|y-x_{k+1}\|^2+\langle \tilde \nabla f(x_{k+1}),y-x_{k+1}\rangle \right\}- \tau_k \delta\\
\geq\min_{y \in Q}\left\{\frac{1}{2}L\|y-x_{k+1}\|^2+\langle \tilde \nabla f(x_{k+1}),y-x_{k+1}\rangle \right\} - \tau_k \delta. \\
\EA\]
By the definition of $y_{k+1}=\tilde T_Q(x_{k+1})$ and using condition (\ref{eq:approx-grad}), we get:
\[\BA{l}
\min_{y \in Q}\left\{\frac{1}{2}L\|y-x_{k+1}\|^2+\langle \tilde \nabla f(x_{k+1}),y-x_{k+1}\rangle \right\} - \tau_k \delta \\
=\frac{1}{2}L\|\tilde T_Q(x_{k+1})-x_{k+1}\|^2+\langle \tilde \nabla f(x_{k+1}),\tilde T_Q(x_{k+1})-x_{k+1}\rangle -  \tau_k \delta\\
\geq\frac{1}{2}L\|y_{k+1}-x_{k+1}\|^2+\langle \nabla f(x_{k+1}),y_{k+1}-x_{k+1}\rangle - 2 \tau_k \delta,\\
\EA\]
and inequality (\ref{eq:lip-ineq}) gives:
\[\BA{l}
\frac{1}{2}L\|y_{k+1}-x_{k+1}\|^2+\langle \nabla f(x_{k+1}),y_{k+1}-x_{k+1}\rangle - 2 \tau_k \delta,\\
\geq f(y_{k+1})-f(x_{k+1})-2\tau_k \delta.\\
\EA\]
Combining these inequalities with the inequality on $\psi_{k+1}$ in (\ref{eq:midineq}), we finally get:
\[
\psi_{k+1}\geq A_{k+1}\left[f(y_{k+1})-(1-\tau_k) g(k,\delta)- 3\tau_{k}\delta\right]
\]
which is the desired result.
\end{proof}

We can use this result to study the convergence of the following algorithm given only approximate gradient information. 
\vskip 1ex
\noindent\line(1,0){370}
\vskip 0ex
\noindent \textbf{Smooth minimization with approximate gradient.}
\vskip 1ex
Starting from $x_0$, the prox center of the set $Q$, we iterate:
\begin{enumerate}
\item compute $\tilde \nabla f(x_k)$,
\item compute $y_k=\tilde T_Q(x_k)$,
\item compute $z_k = \argmin_{x \in Q}\left\{\frac{L}{\sigma}d(x)+\sum_{i=0}^k\alpha_i[ f(x_i)+\langle \tilde \nabla f(x_i),x-x_i\rangle]\right\}$,
\item update $x$ using $x_{k+1}=\tau_k z_k + (1-\tau_k)y_k$,
\end{enumerate}
\line(1,0){370}
\vskip 1ex
Again, because solving for $y_k$ and $z_k$ can often be done very efficiently, the dominant numerical step in this algorithm is the evaluation of $\tilde \nabla f(x_k)$. If the step size sequence $\alpha_k$ satisfies the conditions of Lemma \ref{th:psik}, we can show the following convergence result:
\begin{theorem}\label{th:conv}
Suppose $\alpha_k$ satisfies equation (\ref{eq:alphak}), with the iterates $x_k$ and $y_k$ defined in (\ref{eq:zk-def}) and (\ref{eq:xyz-update}), then for any $k\geq 0$ we have:
\[
f(y_k)-f(x^{\star})\leq \frac{Ld(x^{\star})}{A_k \sigma}+3\delta
\]
where $x^{\star}$ is an optimal solution to problem (\ref{eq:smooth-min}).
\end{theorem}
\begin{proof}
If $\alpha_k$ satisfies the hypotheses of lemma \ref{th:psik} we have:
\[
A_k f(y_k)\leq\psi_k+ A_k g(k,\delta)
\]
where $A_k=\sum_{i=0}^k \alpha_i$ and $g(k,\delta)\leq 3\delta$. Now, because $f(x)$ is convex, we also have:
\[
\psi_k \leq \frac{L}{\sigma}d(x^{\star})+A_k f(x^{\star})+A_k3\delta
\]
which yields the desired result.
\end{proof}

When $d(x^{\star})<+\infty$ (e.g. if ${Q}$ is bounded), if we set the step sequence as $\alpha_k=(k+1)/2$ and $\delta$ to some fraction of the target precision $\epsilon$ (here $\epsilon/6$), $A_k$ grows as $O(k^2)$ and Theorem \ref{th:conv} ensures that the algorithm will converge to an $\epsilon$ solution in less than:
\BEQ \label{eq:maxiter}
\sqrt{\frac{8Ld(x^\star)}{\sigma\epsilon}}
\EEQ
iterations. In practice of course, $d(x^\star)$ needs to be bounded a priori and $L$ is often hard to evaluate. A notable exception is when $f(x)$ is a smooth approximation (as in \cite{Nest03,Nest04a} for example), in which case $L$ is know explicitly as a function of the precision. We have implicitly assumed, as in \cite{Nest83}, that the set $Q$ is simple enough so that the complexity of solving the two minimization subproblems in steps 2 and 3 of the algorithm is low relative to that of approximating the gradient. We also implicitly assumed that the set $Q$ is simple enough so that condition (\ref{eq:approx-grad}) can be checked efficiently. In the numerical experiments of Section \ref{s:numer} for example, steps 2 and 3 are Euclidean projections on the unit box and condition (\ref{eq:approx-grad}) is a simple inequality on the leading eigenvalues of the current iterate.

\section{Semidefinite optimization}
\label{s:semidef} Here, we describe in detail how the results of the previous section can be applied to semidefinite optimization. We consider the following maximum eigenvalue problem:
\BEQ
\label{eq:min-maxeig}\BA{ll}
\mbox{minimize} & \lambda^\mathrm{max}(A^Ty+c)-b^Ty\\
\mbox{subject to} & y \in Q,
\EA\EEQ
in the variable $y\in \reals^m$, with parameters $A\in\reals^{m \times n^2}$, $b\in\reals^m$ and $c\in \reals^{n^2}$. Let us remark that when $Q$ is equal to $\reals^m$, the dual of this program is a semidefinite program with constant trace written:
\BEQ
\label{eq:max-dualsdp}\BA{ll}
\mbox{maximize} & c^Tx\\
\mbox{subject to} & Ax=b\\
& \Tr(x)=1\\
& x \succeq 0,
\EA
\EEQ
in the variable $x\in \reals^{n^2}$, where $\Tr(x)=1$ means that the matrix obtained by reshaping the vector $x$ has trace equal to one and $x\succeq 0$ means that this same matrix is symmetric, positive semidefinite.
\subsection{Smoothing technique}
As in \cite{Nemi04}, \cite{Nest04a}, \cite{dasp04a} or \cite{BenT04} we can find a uniform $\epsilon$-approximation to $\lambda^\mathrm{max}(X)$ with Lipschitz continuous gradient. Let $\mu>0$ and $X\in \symm_n$, we define:
\[
f_\mu(X)=\mu \log \left(\sum_{i=1}^n e^{\lambda_i(X)/\mu}\right)=\mu \log\left( e^\frac{\lambda^{\mathrm{max}}(X)}{\mu}\left(1+\sum_{i=2}^n e^\frac{\lambda_i(X)-\lambda^{\mathrm{max}}(X)}{\mu}\right)\right)
\]
where $\lambda_i(X)$ is the $i^\mathrm{th}$ eigenvalue of $X$. This is also:
\BEQ\label{eq:fmu}
f_\mu(X)=\lambda^{\mathrm{max}}(X) + \mu \log \Tr\left(\exp\left(\frac{X-\lambda^{\mathrm{max}}(X)\mathbf{I}}{\mu}\right)\right)
\EEQ
which requires computing a matrix exponential at a numerical cost of $O(n^3)$. We then have:
\[
\lambda^\mathrm{max}(X) \leq f_\mu(X) \leq \lambda^\mathrm{max}(X) + \mu \log n,
\]
so if we set $\mu=\epsilon/\log n$, $f_\mu(X)$ becomes a uniform $\epsilon$-approximation of $\lambda^\mathrm{max}(X)$. In \cite{Nest04a} it was shown that $f_\mu(X)$ has a Lipschitz continuous gradient with constant:
\[
L=\frac{1}{\mu}=\frac{\log n}{\epsilon}.
\]
The gradient $\nabla f_\mu(X)$ can also be computed explicitly as:
\BEQ\label{eq:grad-fmu}
\frac{\exp\left(\frac{X-\lambda^{\mathrm{max}}(X)\mathbf{I}}{\mu}\right)}{\Tr\left(\exp\left(\frac{X-\lambda^{\mathrm{max}}(X)\mathbf{I}}{\mu}\right)\right)}
\EEQ
using the same matrix exponential as in (\ref{eq:fmu}). Let $\|y\|$ be some norm on $\reals^m$ and $d(x)$ a strongly convex prox-function with parameter $\sigma>0$. As in \cite{Nest04a}, we define:
\[
\|A\|=\max_{\|h\|=1}\|A^Th\|_2,
\]
where $\|A^Th\|_2=\max_i|\lambda_i(A^Th)|$. The algorithm detailed in \cite{Nest03}, where \emph{exact} function values and gradients are computed, will find an $\epsilon$ solution to (\ref{eq:min-maxeig}) after at most:
\BEQ\label{eq:maxiter-sdp}
\frac{2\|A\|}{\epsilon}\sqrt{\frac{\log n}{\sigma}d(y^{\star})}
\EEQ
iterations, each iteration requiring a matrix exponential computation.

\subsection{Spectrum \& expected performance gains}\label{s:expect-perf}
The complexity estimate above is valid when the matrix exponential in (\ref{eq:fmu}) is computed exactly, at a cost of $O(n^3)$. As we will see below, only a few leading eigenvalues are sometimes required to satisfy conditions (\ref{eq:approx-grad}) and obtain a comparable complexity estimate at a much lower numerical cost. To illustrate the potential complexity gains, let us pick a matrix $X\in\symm_n$ whose coefficients are centered independent normal variables with second moment given by $\sigma^2/n$. From Wigner's semicircle law, $\lambda^{\mathrm{max}}(X)\sim 2\sigma$ as $n$ goes to infinity and the eigenvalues of $X$ are asymptotically distributed according to the density:
\[
p(x)=\frac{1}{2\pi \sigma^2}\sqrt{4\sigma^2-x^2},
\]
which means that, in the limit, the proportion of eigenvalues required to reach a precision of $\gamma$ in the exponential is given by:
\[
P_\lambda \triangleq P\left(e^{\frac{\lambda_i(X)-\lambda^{\mathrm{max}}(X)}{\mu}}\leq \gamma\right)=\int_{-2\sigma}^{2\sigma+ \epsilon \frac{\log \gamma}{\log n}}\frac{1}{2\pi \sigma^2}\sqrt{4\sigma^2-x^2}dx.
\]
Since the problems under consideration are relaxations of sparse PCA, we can also consider the case where $X\in\symm_n$ is sampled from the Wishart distribution. In that case, the eigenvalues are distributed according to the Mar\u{c}enko-Pastur distribution (see \cite{Maru67}) and the above proportion becomes:
\[
P_\lambda = P\left(e^{\frac{\lambda_i(X)-\lambda^{\mathrm{max}}(X)}{\mu}}\leq \gamma\right)=\int_{-2\sigma}^{2\sigma+ \epsilon \frac{\log \gamma}{\log n}}\frac{\sqrt{x(4\sigma-x)}}{2\pi x}dx.
\]
With $n=5000$, $\gamma=10^{-6}$ and $\epsilon=10^{-2}$, we get $nP_{\lambda}=2.3$, so the approximations above would suggest that, in theory, it is only necessary to compute about three eigenvalues per iteration to get an approximation with precision $\gamma=10^{-6}$. In practice however, the results of Section \ref{s:numer} show that these rough estimates should be significantly tempered.

\subsection{Global complexity bound}
Let us now focus on the following program:
\BEQ\label{eq:max-eig}
\BA{ll}
\mbox{minimize} & \lambda^\mathrm{max}(A^Ty+c)\\
\mbox{subject to} & \|y\|\leq \beta,\\
\EA\EEQ
where the set ${Q}$ is here explicitly given by :
\[
{Q}=\left\{y\in\reals^p:~\|y\|\leq\beta\right\},
\]
for some $\beta>0$ with $\|.\|$ the Euclidean norm here. We can pick $\|x\|^2/2$ as a prox function for ${Q}$, which is strongly convex with convexity parameter 1. Let $\lambda(X)\in\reals^n$ be the eigenvalues of the matrix $X=A^Ty+c$, in decreasing order, with $u_i(X)\in\reals^n$ an orthonormal set of eigenvectors. The gradient matrix of $\exp(X/\mu)$ is written:
\[
\nabla f_\mu(X)=\left(\sum_{i=1}^n e^{\frac{\lambda_i(X)}{\mu}}\right)^{-1}\sum_{i=1}^n e^{\frac{\lambda_i(X)}{\mu}}u_i(X)u_i(X)^T,
\]
Suppose we only compute the first $m$ eigenvalues and use them to approximate this gradient by:
\[
\tilde \nabla f_\mu(X)=\left(\sum_{i=1}^m e^{\frac{\lambda_i(X)}{\mu}}\right)^{-1}\sum_{i=1}^m e^{\frac{\lambda_i(X)}{\mu}}u_i(X)u_i(X)^T,
\]
we get the following bound on the error:
\[
\|\nabla f_\mu(X)-\tilde \nabla f_\mu(X)\|\leq \frac{\sqrt{2}(n-m) e^{\frac{\lambda_m(X)-\lambda_1(X)}{\mu}}}{\left(\sum_{i=1}^m e^{\frac{\lambda_i(X)-\lambda_1(X)}{\mu}}\right)}.
\]
In this case, with $X=A^Ty-c$ here, condition (\ref{eq:approx-grad}) means that we only need to compute $m$ eigenvalues with $m$ such that:
\BEQ\label{eq:cond-maxeig}
\frac{\sqrt{2}(n-m) e^{\frac{\lambda_m(X)-\lambda_1(X)}{\mu}}}{\left(\sum_{i=1}^m e^{\frac{\lambda_i(X)-\lambda_1(X)}{\mu}}\right)}\leq \frac{\delta}{\sigma^{\mathrm{max}}(A)},
\EEQ
where $\sigma^{\mathrm{max}}(A)$ is the largest singular value of the matrix $A$. Using the result in \cite{Nest04a}, if we define $\|A\|=\max_{\|h\|=1} \|Ah\|_2$ and set $\delta=\epsilon/6$, the algorithm in Section \ref{s:algo} will then converge to an $\epsilon$-solution of problem (\ref{eq:max-eig}) in at most:
\BEQ \label{glob-comp}
\frac{4\|A\|\beta}{\epsilon}\sqrt{\log n}
\EEQ
iterations. This bound on the number of iterations is independent of $m$ in condition (\ref{eq:cond-maxeig}), i.e. the number of eigenvalues required at each iteration. The cost per iteration however varies with problem structure as each iteration requires computing $m$ leading eigenvalues, which can be performed in $O(mn^2)$ operations. Note that partial eigenvalue decompositions only access the matrix through matrix-vector products (see \cite{Leho98}), hence can handle sparse problems very efficiently. The threshold $\delta$ can be adjusted empirically to tradeoff between the number of iterations and the numerical cost of each iteration. Unfortunately, we can't directly infer a bound on $m$ from the structure of $A$, so in the next section we study the link between $m$ and the matrix spectrum in numerical examples.

\section{Examples \& numerical performance} \label{s:numer}
In this section, we illustrate the behavior of the approximate gradient algorithm on various semidefinite optimization problems. Overall, while there appears to be a direct link between problem structure and complexity (i.e. the number of eigenvalues required in the gradient approximation) in the first sparse PCA example discussed below, we will observe on random maximum eigenvalue minimization problems that predicting complexity based on overall problem structure remains an open numerical question in general.

\subsection{Sparse principal component analysis}
\label{ss:spca} Based on the results in \cite{dasp04a}, the problem of finding a sparse leading eigenvector of a matrix $C\in\symm_n$ can be written:
\BEQ\label{eq:variat-prog}
\BA{ll}
\mbox{maximize} & x^TCx\\
\mbox{subject to} & \|x\|_2=1\\
& \Card(x) \leq k,
\EA
\EEQ
where $\Card(x)$ is the number of nonzero coefficients in $x$, and admits the following semidefinite relaxation:
\BEQ
\label{eq:penalized-relax}
\BA{ll}
\mbox{maximize} & \Tr(CX) - \rho \ones^T |X| \ones\\
\mbox{subject to} & \Tr(X)=1\\
& X \succeq 0,
\EA \EEQ
which is a semidefinite program in the variable $X\in \symm^n$,
where $\rho>0$ is the penalty controlling the sparsity of the solution. Its dual is given by:
\BEQ
\label{eq:dual-robust}
\BA{ll}
\mbox{minimize} & \lambda^{\mathrm{max}}(C+U)\\
\mbox{subject to} & |U_{ij}|\leq \rho,\quad i,j=1,\ldots,n,
\EA \EEQ
which is of the form (\ref{eq:min-maxeig}) with
\[
Q=\left\{U\in\symm_n:~ |U_{ij}|\leq \rho,~i,j=1,\ldots,n \right\}.
\]
The smooth algorithm detailed in Section \ref{s:algo} is explicitly described for this problem in \cite{dasp04a} and implemented in a numerical package called DSPCA which we have used in the examples here. To test its performance, we generate a matrix $M$ with uniformly distributed coefficients in $[0,1]$. We let $e\in\reals^{250}$ be a sparse vector with:
\[
e=(1,0,1,0,1,0,1,0,1,0,0,0,\ldots).
\]
We then form a test matrix $C=M^TM+v ee^T$, where $v$ is a signal-to-noise ratio. 

\begin{figure}[htbp]
\begin{center}
\begin{tabular}{cc}
\psfrag{cpu}[t][b]{CPU time (in sec.)}
\psfrag{gap}[b][t]{Duality gap}
\includegraphics[width=0.45 \textwidth]{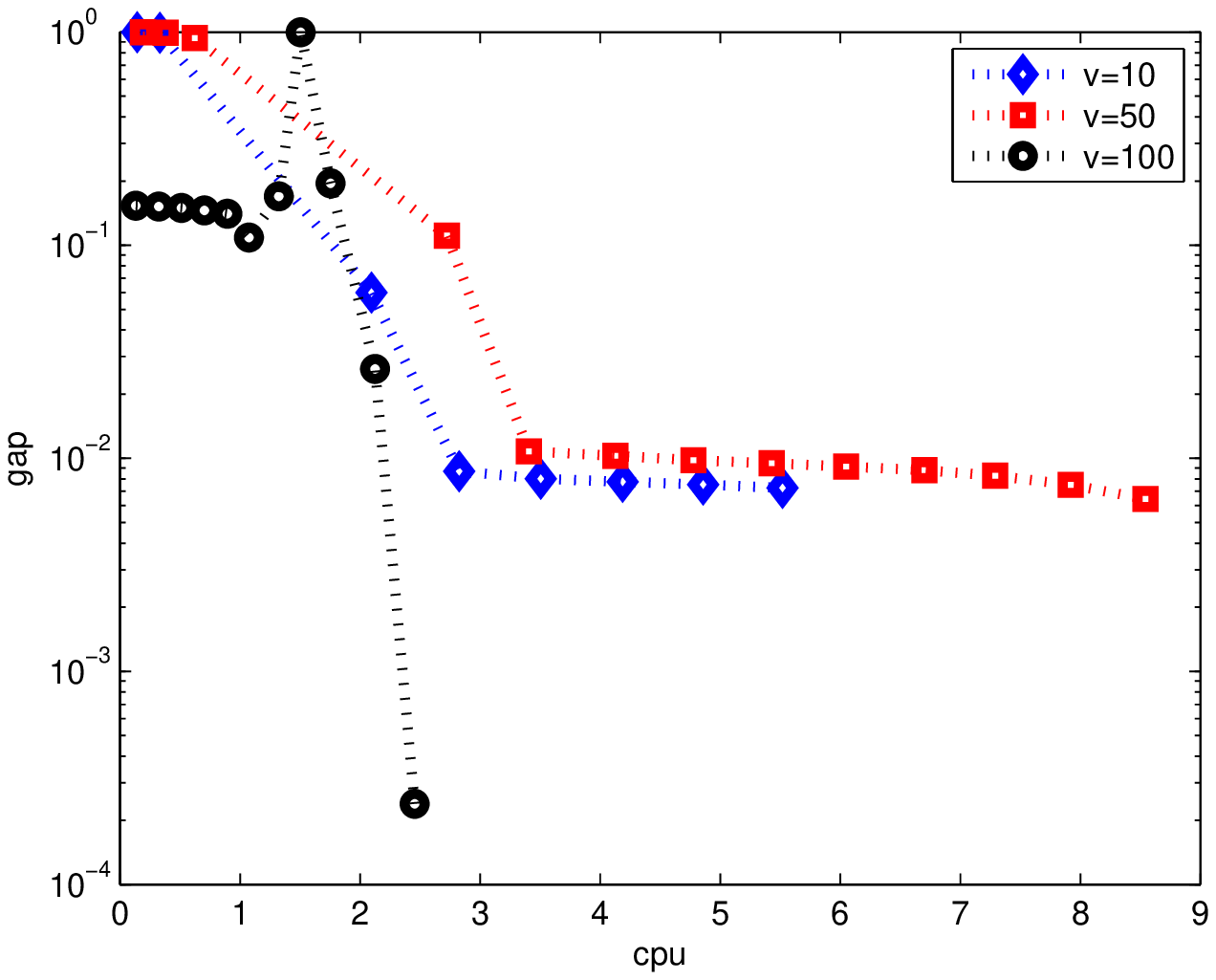}&
\psfrag{time}[t][b]{CPU time (in sec.)}
\psfrag{eigs}[b][t]{\% eigs.}
\includegraphics[width=0.45\textwidth]{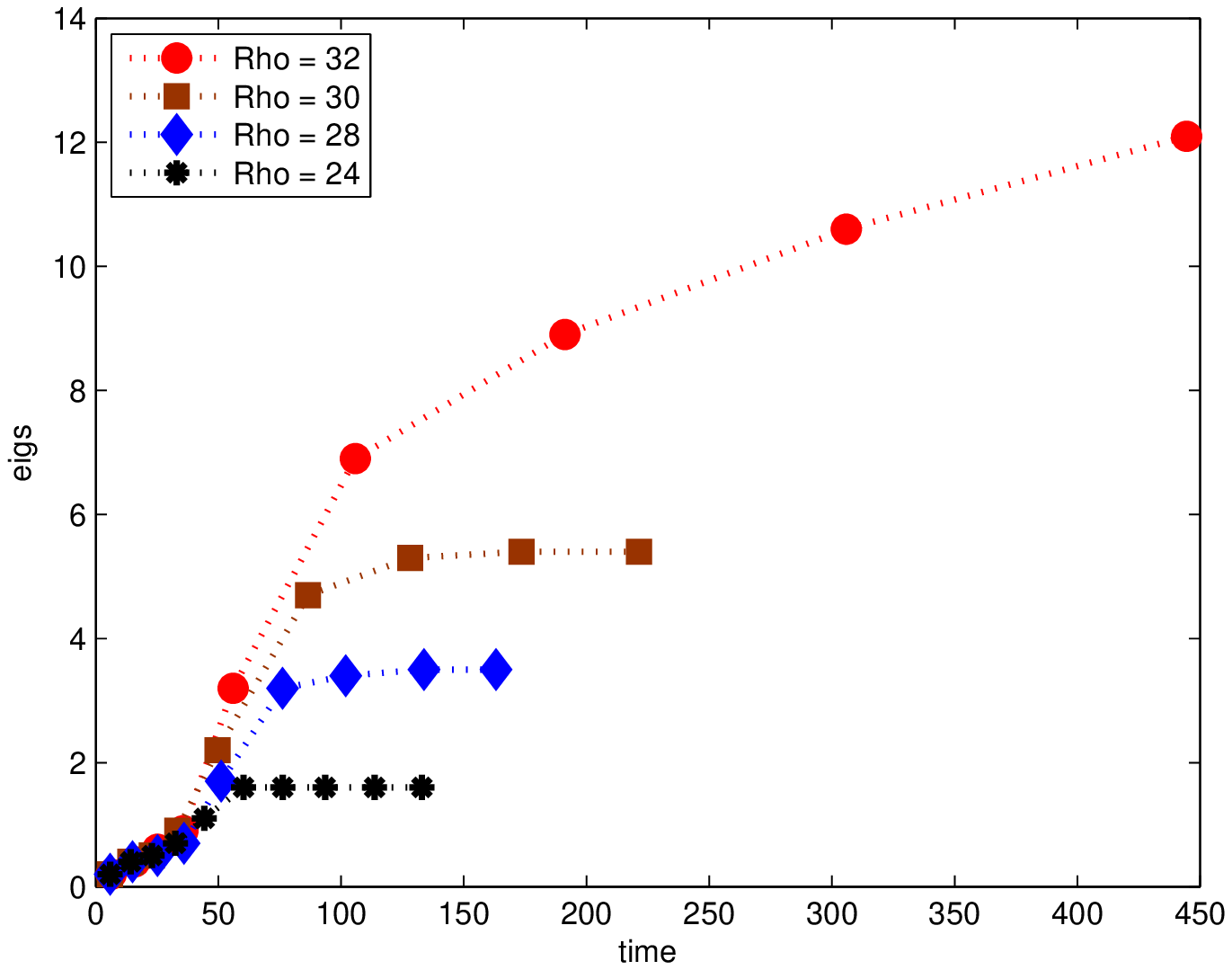}
\end{tabular}
\caption{\textit{Left:} Duality gap versus CPU time for various values of the signal to noise ratio~$v$. \textit{Right:} Percentage of eigenvalues required versus CPU time, for various values of the penalty parameter $\rho$ controlling sparsity. \label{fig:cpu-pca}}
\end{center}
\end{figure}

\begin{figure}[htbp]
\begin{center}
\begin{tabular}{cc}
\psfrag{eps}[t][b]{Precision target}
\psfrag{cpu}[b][t]{CPU time (in sec.)}
\includegraphics[width=0.45 \textwidth]{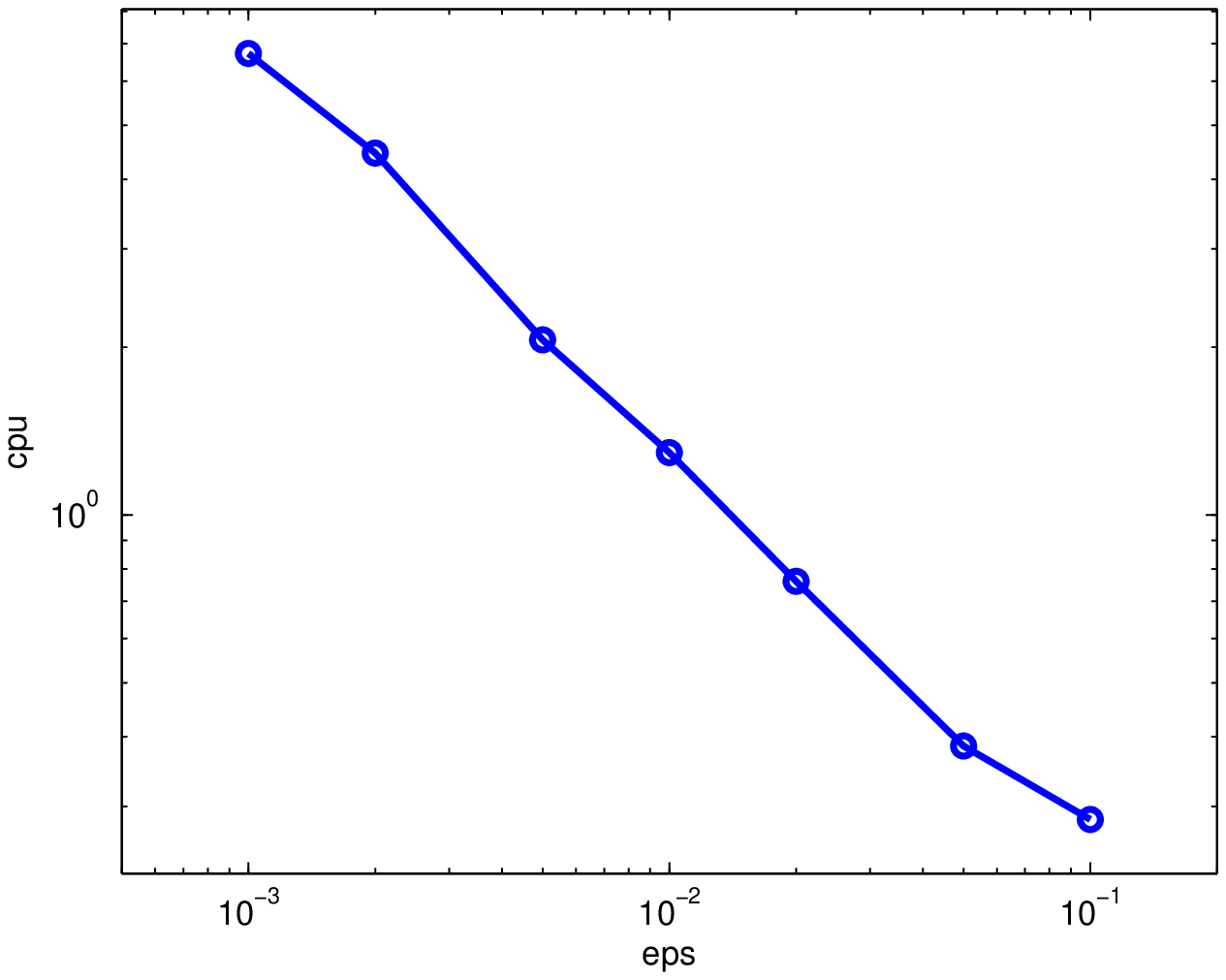}&
\psfrag{eps}[t][b]{Precision target}
\psfrag{perceigs}[b][t]{\% eigs.}
\includegraphics[width=0.45\textwidth]{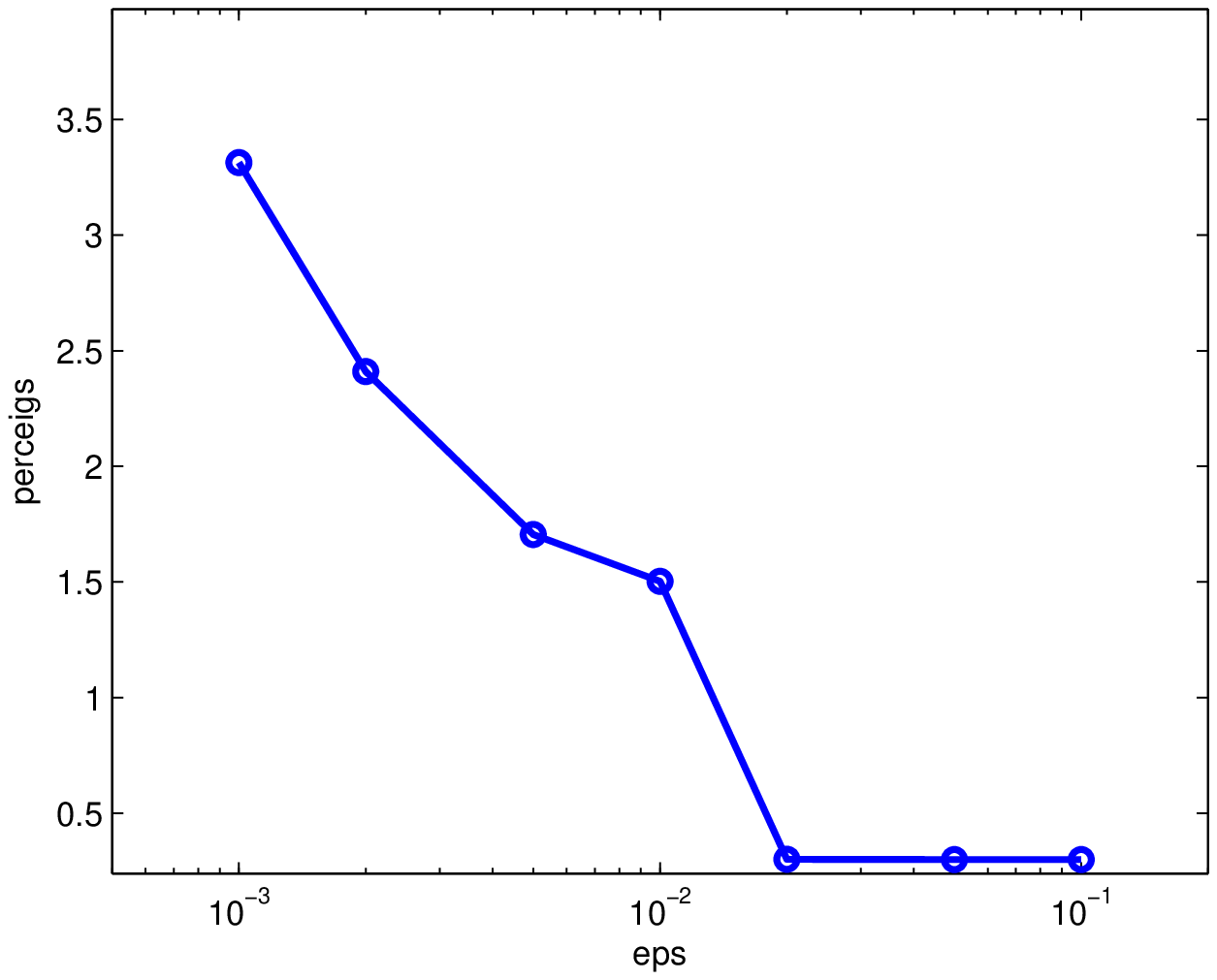}
\end{tabular}
\caption{\textit{Left:} CPU time (in seconds) versus target precision in loglog scale. \textit{Right:} average percentage of eigenvalues required at each iteration versus target precision, in semilog scale. \label{fig:eig-eps}}
\end{center}
\end{figure}

In Figure \ref{fig:cpu-pca} on the left, we plot duality gap versus CPU time used for values of the signal to noise ratio $v$ ranging from 10 to 100. In Figure \ref{fig:cpu-pca} on the right, we plot number of eigenvalues required against computing time using a covariance matrix of dimension $n=500$ sampled from the colon cancer data set in \cite{Alon99}, and a noisy rank one matrix. Finally, we measure total computing time versus problem dimension $n$ on this same data set, by solving problem (\ref{eq:penalized-relax}) for increasingly large submatrices of the original covariance matrix. In each of these examples, we stop after the duality gap has been reduced by $10^{-2}$, which is enough here to identify sparse principal components. In Figure \ref{fig:eig-eps} on the left, we plot computing time versus target precision in loglog scale, on a sparse PCA problem of size 200 extracted from the colon cancer data set. In the previous section, we have seen that precision impacts computing time both through the total number of iterations in (\ref{eq:maxiter-sdp}) and through condition (\ref{eq:cond-maxeig}) on the number of eigenvalues required in the gradient approximation. In this example, we observe that CPU time increases a little bit slower than the upper bound of $O(1/\epsilon)$  given in (\ref{fig:eig-eps}). In Figure \ref{fig:eig-eps} on the right, we plot the average percentage of eigenvalues required at each iteration versus target precision, in semilog scale. We observe, on this example of dimension 200, that for low target precisions, one eigenvalue is often enough to approximate the gradient, but that this number quickly increases for higher precision targets. Note that in all cases, the precision targets are significantly lower than those achieved by interior point methods (usually at least $10^{-8}$) but the cost per iteration and storage requirements of the first-order algorithms detailed here are also significantly lower.

In Table \ref{tab:cputime}, we then compare total CPU time using a full precision matrix exponential, against CPU time using only a partial eigenvalue decomposition to approximate this exponential. Note that other classic methods for computing the matrix exponential such as Pad\'e approximations (see \cite{Mole03}), did not provide a significant performance benefit and are not included here. Both exact and approximate gradient codes are fully written in C, with partial eigenvalue decompositions computed using the FORTRAN package ARPACK (see \cite{Leho98}) with calls to vendor optimized BLAS and LAPACK for matrix operations. To improve stability, the size of the Lanczos basis in ARPACK was set at four times the number of eigenvectors required. We observe that, on these problems, the partial eigenvalue decomposition method is about ten times faster.

\begin{table}[ht]
\begin{center}
\begin{tabular}{r|ccc}
& $n=100$ & $n=200$ & $n=500$ \\
\hline
Rank one, Full & 3.2 & 8.0 & 14.7 \\
Rank one, \textbf{Partial} & 0.4 & 0.75 & 1.6 \\
Colon, Full & 2.6 & 18.1 & 274.3 \\
Colon, \textbf{Partial} & 0.3 & 1.3 & 17.7 \\
\end{tabular}
\end{center}
\caption{CPU time (in seconds) versus problem dimension $n$ for full and partial eigenvalue matrix exponential computations.\label{tab:cputime}}
\end{table}

\subsection{Matrix structure and complexity: open numerical issues}
The previous section showed how the spectrum of the current iterate impacts the complexity of the algorithm detailed in this paper: a steeply decreasing spectrum allows fewer eigenvalues to be computed in the matrix exponential approximation, and a wider gap between eigenvalues improves the convergence rate of these eigenvalue computations. In this section, we study the number of eigenvalues required in randomly generated maximum eigenvalue minimization problems. Because of the measure concentration phenomenon, there is nothing really random about the spectrum of large-scale, naively generated semidefinite optimization problems so we begin by detailing a simple method for generating random matrices with a given spectrum. 

\paragraph{Generating random matrices with a given spectrum}
Suppose $X\in\symm_n$ is a matrix with normally distributed coefficients, $X_{ij}\sim\mathcal{N}(0,1)$, $i,j=1,\ldots,n$. If we write its QR decomposition, $X=QR$ with $Q,~R\in \reals^{n \times n}$, then the orthogonal matrix $Q$ is Haar distributed on the orthogonal group $\mathcal{O}_n$ (see \cite{Diac03} for example). This means that to generate a random matrix with given spectrum $\lambda\in\reals^n$, we generate a normally distributed matrix $X$, compute its QR decomposition and the matrix $Q\diag(\lambda)Q^T$ will be uniformly distributed on the set of symmetric matrices with spectrum $\lambda$.

\begin{figure}[htbp]
\begin{center}
\psfrag{gap}[c][c]{\small{gap}}
\psfrag{perc}[b][b]{\small{\% eig.}}
\psfrag{Gaussian}[c][t]{\bf Gaussian}
\psfrag{wishart}[c][t]{\bf Wishart}
\psfrag{uniform}[c][t]{\bf Uniform}
\psfrag{uniformrkone}[c][t]{\bf Uniform +1}
\includegraphics[width= \textwidth]{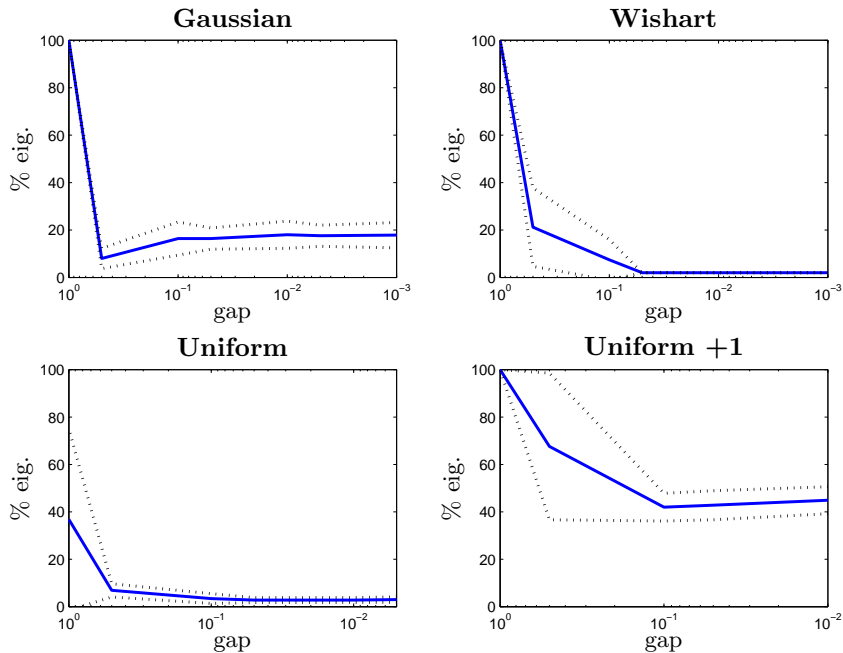}
\caption{Average percentage of eigenvalues required (solid line) versus duality gap on randomly generated maximum eigenvalue minimization problems, for various problem matrix distributions. Dashed lines at plus and minus one standard deviation.
\label{fig:prop-maxeig}}
\end{center}
\end{figure}

\paragraph{Maximum eigenvalue minimization} We now form random maximum eigenvalue minimization problems, then study how the number of required eigenvalues in the gradient computation evolves as the solution approaches optimality. We solve the following problem:
\[\BA{ll}
\mbox{minimize} & \lambda^\mathrm{max}(A^Ty+c)\\
\mbox{subject to} & \|y\|\leq \beta,\\
\EA\]
in the variable $y \in \reals^m$, where $c\in\reals^{n^2}$, $A\in\reals^{m \times n^2}$ and $\beta>0$ is an upper bound on the norm of the solution.
In Figure \ref{fig:prop-maxeig} we plot percentage of eigenvalues required in the gradient computation versus duality gap for randomly generated problem instances where $n=50$ and $m=25$. The first two plots use data matrices with Gaussian and Wishart distributions, whose spectrum are distributed according to Wigner's semicircle law and the Mar\u{c}enko-Pastur distribution respectively. The last two plots use the procedure described above to generate matrices with uniform spectrum on $[0,1]$, and uniform spectrum on $[0,1]$ with one eigenvalue set to 5. We observe that the number of eigenvalues required in the algorithm varies significantly with matrix spectrum.

\paragraph{Problem structure and effective complexity} 
The results on sparse PCA in $\S\ref{ss:spca}$ and on the random problems of this section show that problem structure has a significant impact on performance. Predicting how many eigenvalues will be required at each iteration based on structural properties of the problem is an important but difficult question. In particular, the number of eigenvalues required in the Gaussian case is much higher than what the asymptotic analysis in Section \ref{s:expect-perf} predicted. Furthermore, in the sparse PCA example, complexity seems to vary with problem structure somewhat intuitively: higher signal to noise ratio means lower complexity and a higher sparsity target means higher complexity. However, this is not the case in the random problems studied here, two unstructured problems (uniform and Wishart) have low complexity while one requires computing many more eigenvalues per iteration (Gaussian) and a more structured example (uniform plus rank one) also requires many eigenvalues. Overall then, predicting effective complexity (i.e. the number of eigenvalues required at each iteration) based on problem structure remains a difficult open question at this point.

Also, it is well known empirically (see \cite{Over92a}, \cite{Lewi96} and \cite{Pata98} among others) that the largest eigenvalues of $A^Ty-c$ in (\ref{eq:min-maxeig}) tend to coalesce near the optimum, thus potentially increasing the number of eigenvalues required when computing $\tilde \nabla f(x)$ and the number of iterations required for computing leading eigenvalues (see \cite{Golu90} for example), but in these references too, no a priori link between coalescence and problem structure is established. This coalescence phenomenon is never apparent in the numerical examples studied here, perhaps because it only appears at the much higher precision targets reached by interior point methods.

\section*{Acknowledgements}
The author would like to thank Noureddine El Karoui for very useful comments and to acknowledge support from NSF grant DMS-0625352, ONR grant number N00014-07-1-0150, a Peek junior faculty fellowship and a gift from Google, Inc.

\bibliographystyle{siam}
\bibliography{Mainperso}
\end{document}